\documentclass[12pt,a4paper,reqno]{amsart}
\usepackage{preemble}
\usepackage{tikz, enumerate}
\usepackage{xcolor}

\date{\today}

\author{Stephen Cantrell}
\address{Mathematical Institute, University of St Andrews, North Haugh, St Andrews, KY16 9SS, Scotland}
\email{sjc33@st-andrews.ac.uk}

\date{\today. 
}

\title[Quantitative spectral rigidity for manifolds and Outer Space]{Quantitative spectral rigidity for negatively curved manifolds and Outer Space}

\begin{document}
\maketitle

\begin{abstract} 
Let $M$ be a closed manifold with fundamental group $G$. A subset of conjugacy classes $E \subset \conj(G)$ is said to be spectrally rigid if when two negatively curved Riemannian metrics have the same marked length spectra on $E$, then their length spectra agree on all conjugacy classes. In this work we study spectrally rigid subsets from a quantitative perspective. 

We show that, if the density of $E \subset \conj(G)$, measured with respect to a fixed reference metric, decays more slowly than $T^{-1/2}$ then $E$ is spectrally rigid. We establish analogous rigidity results for  Culler-Vogtmann Outer Space and show that the $T^{-1/2}$ threshold is sharp in this setting.
\end{abstract}

\section{Introduction}
Let $M$ be a closed manifold equipped with a smooth Riemannian metric $\gfr$ of negative sectional curvature. Each non-trivial conjugacy class $[g] \in \conj(G)$ in the fundamental group $G = \pi_1(M)$ corresponds to a unique closed geodesic on $(M,\gfr)$ and the marked length spectrum $\ell_\gfr : \conj(G) \to \R$ assigns to each conjugacy class the length of its associated closed geodesic. The marked length spectrum rigidity conjecture (due to Burns and Katok \cite{burnskatok}) predicts that the isometry class of $(M,\gfr)$ is uniquely determined, amongst all negatively curved Riemannian metrics on $M$, by its marked length spectrum. 
This conjecture was established for surfaces independently by Otal \cite{Otal} and Croke \cite{Croke}, while in higher dimensions partial results are known; see for example \cite{hamenstadt, g.l}.

In recent years there has been a growing interest in identifying spectrally rigid sets: collections of conjugacy classes whose lengths determine the full marked length spectrum. More precisely, if $(M,\gfr)$ is a closed negatively curved Riemannian manifold with fundamental group $G$, a subset $E \subset \conj(G)$ is called spectrally rigid if when two negatively curved metrics on $M$ have marked length spectra that agree on $E$  then they agree on all $\conj(G)$.

A number of spectrally rigid sets with strong geometric or algebraic structure are now known. Katok \cite{katok} showed that, for surfaces, closed geodesics in a fixed homology class are spectrally rigid for metrics within a fixed conformal class; this restriction was later removed and the result extended to higher dimensions by Gogolev and Rodriguez Hertz \cite{gogolev.hertz}. The author and Reyes established spectral rigidity for closed geodesics lying in a coamenable subgroup \cite{cantrell.reyes.lsr}, while analogous results for Anosov flows appear in \cite{sharp.liv, sharp.pollicott.liv}. Hao \cite{yanlong} showed that subgroups with full limit set in the universal cover are spectrally rigid.

Beyond such structured examples, it is natural to ask how `large' a set must be in order to guarantee its spectral rigidity. This quantitative perspective was first studied in \cite{noelle} and \cite{dilsavor.reber} where size was measured by asymptotic growth with respect to geometric length. To describe these results, fix a negatively curved metric $\gfr$ on $M$, and for a subset $E \subset \conj(G)$, define its density by
\[
E_\gfr(T)=\frac{\#\{[g]\in E:\ell_\gfr[g]<T\}}{\#\{[g]\in\mathrm{conj}(G):\ell_\gfr[g]<T\}}
\]
for $T>0$.

In \cite{noelle}, it was shown that for surfaces, if the complement $E^c$ grows subexponentially then $E$ is spectrally rigid. This was sharpened in \cite{dilsavor.reber}  where it was  proved that if $E$ has positive upper density  (i.e. $\limsup_{T\to\infty} E_\gfr(T) > 0$) then $E$ is spectrally rigid.   On the other hand, it was shown in \cite{cantrell.lsr} and \cite{yanlong.3} that for any prescribed decay rate,  there exists an infinite spectrally rigid sets $E$ whose density $E_\gfr(T)$ decays faster than this rate. We also note that spectral rigidity is not forced below full exponential growth rate in the following sense: for any $\alpha >0$ there exists a non-spectrally rigid set such that $\limsup_{T\to\infty} e^{\alpha T} E_\gfr(T) = \infty$, see \cite[Section 1]{cantrell.lsr}. These results naturally lead us to the following question:\\

\textit{How large must a set of conjugacy classes be in order to guarantee its spectral rigidity? In particular, is rigidity forced by a growth condition that is weaker than positive (upper) asymptotic density?} \\

Our first main result answers this question affirmatively. 

\begin{theorem}\label{thm.manifolds}
Let $(M,\mathfrak{g})$ be a smooth, closed, negatively curved Riemannian manifold with fundamental group $G$. If a set $E \subset \conj(G)$  satisfies
\[
\limsup_{T\to\infty} \sqrt{T} \, E_\gfr(T) = \infty
\]
then it is spectrally rigid.
\end{theorem}

In fact, our proof yields a stronger conclusion: for $E$ satisfying the growth condition stated in theorem and for any two negatively curved metrics $\gfr_1, \gfr_2$ on $M$, the difference $\ell_{\gfr_1} - \ell_{\gfr_2}$ is unbounded on $E$. We expect that rigidity should hold under weaker growth assumptions. Indeed, it seems reasonable that if there exists a polynomial $P$ such that  $\limsup_{T\to\infty} P(T) E_\gfr(T) = \infty$, then $E$ should be spectral rigid. Obtaining such improvements appears to rely on understanding subtle arithmetic properties of length spectra; specifically, whether the subgroup generated  $\{\ell_{\gfr_1}[g] - \ell_{\gfr_2}[g] : g \in G\}$ has finite rank or not (see Lemma \ref{lem.nonarith} where we establish a result of this flavour under additional assumptions). Such problems are notoriously challenging and are closely related to mixing rates of the geodesic flow \cite{dolgopyat}. 
\begin{remark}
It is natural to ask for a dynamical version of Theorem \ref{thm.manifolds} in the spirit of \cite[Theorem 1.1]{dilsavor.reber}. That is, can we obtain a Livsics Theorem for Anosov flows where we assume that two potentials have the same periodic orbit data on a sufficiently large (but decaying in terms of density) collection of periodic orbits? We believe that this should be possible for Anosov flows that are exponentially mixing but we will not pursue this result in the current paper.
\end{remark}
After establishing Theorem \ref{thm.manifolds}, we turn to the analogous problem in Culler-Vogtmann Outer Space. 
Let $F_N$ be the free group of rank $N \ge 2$. Outer Space $\cv_N$ consists of free, minimal, discrete, isometric actions of $F_N$ on $\R$-trees, considered up to equivariant isometry. Each $T \in \cv_N$ (where $T$ denotes the $\R$-tree corresponding to the action) is determined by its translation length function
\[
\|g\|_T=\inf_{x\in T}d_T(x,gx),
\]
where $d_T$ denotes the metric on $T$ \cite{cullermorgan}.
As in the manifold setting, one may ask which subsets of conjugacy classes determine the full translation length function. Finite spectrally rigid sets do not exist for $\cv_N$ \cite{smillie.vogtmann, cls}, but various infinite examples are known. Coamenable subgroups and conjugacy classes in fixed homology classes are spectrally rigid \cite{cantrell.reyes.lsr}, as well as orbits of sufficiently rich collections of automorphisms \cite{many}. On the other hand, finitely generated infinite-index subgroups  \cite{ray.1} and the orbits of individual automorphisms \cite{ray.2} are not spectrally rigid. 
An interesting result of Kapovich shows that typical geodesic rays in Cayley graphs are spectrally rigid \cite{Kap}.

Fix a finite generating set $S$ for $F_N$. For a conjugacy class $[g] \in \conj(F_N)$, let $\ell_S[g]$ denote its stable word length (see Section \ref{sec.ghyp}), and define the density by
\[
E_S(T)=\frac{\#\{[g]\in E:\ell_S[g]<T\}}{\#\{[g]\in\mathrm{conj}(F_N):\ell_S[g]<T\}}
\]
for $T>0$.
Our main result for Outer Space gives a sharp threshold for spectral rigidity in terms of the decay of $E_S(T)$. We have chosen to count with respect to $S$ as this is in keeping with the search for `small' spectrally rigid sets initiated by Kapovich in \cite{Kap}.
\begin{theorem} \label{thm.outermain}
Let $F_N$ be the free group with rank $N\ge 2$ equipped with a finite generating set $S$. If $E \subset \conj(F_N)$ satisfies
\[
\limsup_{T\to\infty} \sqrt{T} \, E_S(T)= \infty, 
\]
then $E$ is spectrally rigid for $\cv_N$. Furthermore, this result is optimal: for each $N \ge 2$ there exist a set $Q \subset \conj(F_N)$ and a generating set $S$ for $F_N$ such that
\[
0 < \limsup_{T\to\infty} \sqrt{T} \, Q_S(T)  < \infty, 
\]
and $Q$ is not spectrally rigid for $\cv_N$.
\end{theorem}

More generally, we deduce these results from a rigidity theorem for cocompact actions on $\text{CAT}(-1)$ metric spaces. As in the manifold case, the proof also shows that the length differences between any two such actions must be unbounded on sets satisfying the growth condition above.


\subsection*{Organisation of the paper.}
We begin by reviewing background material from geometric group theory and symbolic dynamics. In Section \ref{sec.outer} we prove our results for Outer Space and more generally $\textnormal{CAT}(-1)$ spaces. The final section is devoted to the  proof of Theorem \ref{thm.manifolds} for negatively curved manifolds.

\subsection*{Acknowledgements}
The author would like to thank Caleb Dilsavor, Natalia Jurga and Richard Sharp for helpful comments on a preliminary draft of this paper. We also thank the referee for their comments and suggestions. We feel that these have greatly improved the presentation of this work.

\subsection*{AI disclosure}
AI was not used for this project, neither to assist with the writing nor to obtain mathematical statements/proofs.

\section{Preliminaries}
In this section we  briefly recall the background material from geometric group theory and symbolic dynamics that will be used throughout the paper. For a more detailed account of the material see \cite{cdp} and \cite{ParryPollicott}.
\subsection{Coarse geometry}
\subsubsection{Outer space}
Let $F_N$ be the free group of rank $N \ge 2$. The Culler-Vogtmann Outer Space $\cv_N$ is defined as follows. Suppose that $F_N$ acts by isometries on an $\R$-tree $T$ (a $0$-hyperbolic metric space, see Section \ref{sec.ghyp} below) and that this action is free, minimal and discrete.  We will abuse notation and denote this action by $T$. 

The space $\cv_N$ consists of equivalence classes of such actions where $T_1\sim T_2$ if there exists an $F_N$ equivariant isometry from $T_1$ to $T_2$.

Given $T \in \cv_N$, let
\[
\|\cdot\|_T : F_N \to \R^+, \, \, \|g\|_T = \inf_{x \in T} d_T(x,gx)
\]
be the translation length function which is constant on conjugacy classes. Given $T_1, T_2 \in \cv_N$ we have that $T_1 = T_2$ if and only if the functions $\|\cdot\|_{T_1}$ and $\|\cdot\|_{T_2}$ agree on all conjugacy classes \cite{cullermorgan}. See \cite{Vogtmannsurvey} for more details on Outer Space.

\subsubsection{General actions on hyperbolic metric spaces}\label{sec.ghyp}
For a pseudometric space $(X, d)$ (i.e. $d$ satisfies the triangle inequality, is symmetric and has the property that $d(x,x)=0$) the Gromov product is defined by
\[
 (g|h)_f =\frac{1}{2}( d(f, g)+d(f, h)-d(g, h)) \ \ \text{for $g, h, f \in X$}.
\]
We say that $(X, d)$ is $\d$-hyperbolic for some $\delta \ge 0$ if
\begin{equation}
(g|h)_f \ge \min\left\{(g|k)_f, (h|k)_f\right\}-\d \quad \text{for all $g, h, f, k \in X$}.
\end{equation}
We call a pseudometric space hyperbolic if it is $\d$-hyperbolic for some $\d \ge 0$.
A hyperbolic group is a finitely generated group $G$ such that, when equipped with the word distance $d_S(g,h) = |g^{-1}h|_S$ associated to a finite generating set $S$, it becomes a hyperbolic metric space. Throughout this work we will always work with groups $G$ that are non-elementary (i.e. not virtually cyclic) and hyperbolic.
 
We say that two pseudometrics $d,d_\ast$ on $G$ are quasi-isometric if there exist constants $L > 0$ and $C \ge 0$ such that
\[
L^{-1} \, d(g,h)-C \le d_\ast(g, h) \le L \, d(g, h)+C \ \text{ for all $g, h \in G$}.
\]
Throughout this work,  $\Dc_G$ will denote the set of pseudometrics which are left-invariant, hyperbolic and quasi-isometric to some (equivalently, any) word metric on $G$. Given $d \in \Dc_G$, we will write $(\cdot|\cdot)_d : G \times G \to \R$ for the Gromov product of $d$ based at $o$, the identity of $G$. We write $(\cdot|\cdot)_S$ when the corresponding metric is the word metric for the finite generating set $S$. If $d_1, d_2 \in \Dc_G$ then there exists $\lambda, C > 1$ such that
\begin{equation}\label{eq.gromprod}
\lambda^{-1} (g|h)_{d_1} - C \le (g|h)_{d_2} \le \lambda (g|h)_{d_1} + C
\end{equation}
for all $g,h \in G$ \cite[Prop. 15 (i), Ch. 5]{cdp}. We note that members of $\Dc_G$ are not necessarily geodesic, but are roughly geodesic: given $d \in \Dc_G$ there exists a constant $\alpha >0$ such that for any $g,h \in G$ there exists a sequence of elements $g_n \in G$ with $g_0 = g, g_n = h$ and $|d(g_i,g_j) - |j-i|| < \alpha$ for any $i,j \in \{0,\ldots, n\}$ \cite[Proposition 5.6]{BonkSchramm}.

Given $d\in\Dc_G$ we define the corresponding translation length function $\ell_d : \conj(G) \to \R$ by
\[
\ell_d[g] = \lim_{n\to\infty} \frac{d(o,g^n)}{n}
\]
where $o$ denotes the identity in $G$. Two pseudometrics $d,d_\ast \in \Dc_G$ have proportional translation length functions if and only if they are roughly similar: there exists $\tau, C >0$ such that $|d(g,h) - \tau d_\ast(g,h)| < C$ for all $g,h \in G$ \cite{FurmanCoarse}, \cite{cantrell.tanaka.1}. When a hyperbolic group $G$ acts properly, cocompactly by isometries on a geodesic, hyperbolic metric space $(X,d_X)$ we obtain a pseudometric in $\Dc_G$ as follows: fix a base point $o \in X$ and define, for $g,h\in G$, $d(g,h) = d_X(go,ho)$. We define the  stable translation function $\ell_X$ associated to this action to be $\ell_d$. We recall that the exponential growth rate $v_d$ of a pseudometric $d \in \Dc_G$ is the exponential growth rate of the size of $T$ balls as $T$ tends to infinity. This agrees with the growth rate of the number of conjugacy class of $\ell_d$ size at most $T$:
\[
v_d := \limsup_{T\to\infty} \frac{1}{T} \log \#\{ g \in G: d(o,x) < T\} = \limsup_{T\to\infty} \frac{1}{T} \log \#\{ [g] \in \conj(G): \ell_d[g] < T\}. 
\]

When studying spectral rigidity for Riemannian metrics we will work with pseudometrics in $\Dc_G$ that satisfy strong regularity properties at infinity (i.e. on the Gromov boundary).
\begin{definition}
A hyperbolic pseudometric $d$ in $G$ is called strongly hyperbolic if there exist positive constants $c, R_0 > 0$ such that for all $R \ge R_0$, and all $g, g', h, h' \in G$,
if
$ d(g, h)-d(g, g')+d(g', h')-d(h, h') \ge R,$
then
\[
|d(g, h)-d(g', h)-d(g, h')+d(g', h')| \le e^{-c R}.
\]
\end{definition}
Recall that hyperbolic groups admits a compactification by its Gromov boundary $\partial G$. This boundary consists of equivalence classes
of divergent sequences. A sequence of group elements $\{g_n\}_{n=0}^\infty$ diverges
if $(g_n|g_m)_o$ diverges as $\min\{n, m\}$ tends to infinity (here $o$ is the identity in $G$ and the Gromov product is with respect to a fixed pseudometric in $\Dc_G$). Two divergent sequences $\{g_n\}_{n=0}^\infty$ and $\{h_n\}_{n=0}^\infty$ are equivalent if 
$(g_n|h_m)_o$ diverges as $\min\{n, m\}$ tends to infinity. 

Let $d \in \Dc_G$ be strongly hyperbolic. The Busemannn function $\beta_d : G \times \partial G \to \R$ associated to $d$ is defined
\begin{equation}\label{eq.limit}
\beta_d(g,\x) = \lim_{n\to\infty} d(g,\x_n) - d(o,\xi_n)  \ \ \text{ where $\xi =  \{ \x_n\}_{n=0}^\infty \in \partial G$. }
\end{equation}
The fact that the above limit exists relies on the fact that $d$ is strongly hyperbolic. Furthermore, for each fixed $g \in G$, the $\beta_d(g,\cdot)$ is H\"older continuous in its second variable: there exists $C, \alpha >0$ such that
\[
|\beta_d(g,\xi) - \beta_d(g,\eta)| \le C e^{-\alpha (\xi|\eta)_o},
\] 
where $(\xi|\eta)_o$ is the natural extension (which we won't define here: see \cite[Section  2.2]{cantrell.tanaka.1}) of the Gromov product of $d$ to the boundary. This follows directly from the definition of strong hyperbolicity.


\subsubsection{Manhattan curves}

To prove our main theorems it will be useful to consider Manhattan curves for pairs of pseudometrics in $\Dc_G$. Fix $d,d_\ast \in \Dc_G$. The Manhattan curve $\theta_{d_\ast/d}$ is defined to be the abscissa of convergence of 
\[
\sum_{x \in G} e^{-td_\ast(o,x) - sd(o,x)}
\]
for fixed $t$ as $s$ varies. That is, $\theta_{d_\ast/d}(t)$ is the unique real number for which the above series converges when $s > \theta_{d_\ast/d}(t)$ and diverges when $s < \theta_{d_\ast/d}(t)$. The Manhattan curve can also be obtained as
\begin{equation}\label{eq.count}
\theta_{d_\ast/d}(t) = \lim_{T\to\infty} \frac{1}{T} \log \sum_{|d(o,x) - T| <R_0} e^{-t d_\ast(o,x)} =  \lim_{T\to\infty} \frac{1}{T} \log \sum_{|\ell_{d}[g] - T| <R_0} e^{-t \ell_{d_\ast}[g]} 
\end{equation}
for any fixed sufficiently large $R_0 >0$ \cite[Corollary 2.10 and Proposition 3.1]{cantrell.tanaka.1}. Here in the first expression for $\theta_{d_\ast/d}$ we count over group elements whereas in the second we count over conjugacy classes. We note that, although the results in \cite{cantrell.tanaka.1} are stated for genuine metrics, the same proofs apply to pseudometrics.

The curve $\theta_{d_\ast/d}$ encodes a surprising amount of information about the pseudometrics $d_\ast,d$. We summaries some of the key properties of $\theta_{d_\ast/d}$  in the following result.

\begin{theorem} [Theorem 1.1 and 1.2 of \cite{cantrell.tanaka.1}]\label{thm.Man}
Fix $d,d_\ast \in \Dc_G$. Then the Manhattan curve $\theta_{d_\ast/d}$ is $C^1$ and the limit
\[
\tau(d_\ast/d) := \lim_{T\to\infty} \frac{1}{\#\{x \in G: d(o,x) < T\}} \sum_{d(o,x)<T} \frac{d_\ast(o,x)}{T}
\]
exists and satisfies $\tau(d_\ast/d) \ge v_d/v_{d_\ast}$. Furthermore, the following are equivalent:
\begin{enumerate}
\item $\theta_{d_\ast/d}$  is a straight line;
\item $d$ and $d_\ast$ are roughly similar; and,
\item $\tau(d_\ast/d) = v_d/v_{d_\ast}$.
\end{enumerate}
\end{theorem}

The following result shows that the Manhattan curve $\theta_{d_\ast/d}$ can be seen as a curve that interpolates between the pseudometrics $d,d_\ast$. 
\begin{proposition}[Proposition 4.1 \cite{cantrell.reyes.man}]\label{lem.bound}
Take $d,d_\ast \in \Dc_G$. Then for any $t \in \R$ there exists $C_t > 0$ such that the function $d_t : G \times G \to \R$ defined by
\[
d_t(g,h) = 
     \begin{cases}
       td_\ast(g,h) + \theta_{d_\ast/d}(t)d(g,h)+ C_t  &\text{ if $g \neq h$}\\
       0 &\text{ otherwise}
     \end{cases}
\]
determines a pseudometric $d_t \in \Dc_G$.
\end{proposition}
We note that, by the definition of the Manhattan curve, the pseudometrics $d_t$ have exponential growth rate $1$.

\subsection{Dynamical systems and probability theory}
\subsubsection{Symbolic dynamics}
We briefly recall some introductory material from symbolic dynamics. See \cite{ParryPollicott} for a comprehensive treatment.
Let $A$ be the $k \times k$ matrix with entries $0$ or $1$. 
We use the notation $A(i,j)$ to denote the $(i,j)$-th entry of $A$. The one-sided subshift of finite type associated to $A$ is the space
\[
\Sigma_A = \{(x_n)_{n=0}^{\infty} : x_n \in \{1,2,...,k\}, A(x_n,x_{n+1})=1, n \in \mathbb{Z}_{\ge 0}\}.
\]
Given $x\in\SS_A$ we write $x_n$ for the $n$-th coordinate of $x$. 
The shift map $\sigma: \Sigma_A \rightarrow \Sigma_A$ sends $x$ to $y= \sigma(x)$ where $y_n=x_{n+1}$ for all $n \in \mathbb{Z}_{\ge 0}$. Given a finite sequence $x_0, \ldots, x_{n-1} \in \{1,\ldots,k\}$ we write $[x_0, \dots, x_{n-1}]$ for the corresponding cylinder set i.e. $[x_0, \dots, x_{n-1}]$ denotes the collection of sequences in $\Sigma_A$ that agree with $(x_0,x_1,\ldots,x_{n-1})$ on the $0$-th to $(n-1)$-st indices. We equip $\Sigma_A$ with the topology generated with by all cylinder sets.

Given a function $f : \Sigma_A \to \R$ we use, for each $n\ge 1$, the notation $f^n(x)$ to denote the $n$-th Birkhoff sum of $f$ evaluated at $x \in \Sigma_A$: $f^n(x) = f(x) + f(\sigma(x)) + \cdots + f(\sigma^{n-1}(x))$.

We write $\Mcc(\SS_A,\sigma)$ for the set of $\sigma$-invariant Borel probability measures on $\SS_A$ and for $\lambda \in \Mcc(\SS_A,\sigma)$ let $h(\s, \lambda)$ denote the measure theoretical entropy of $(\SS_A, \s, \lambda)$ (see Section 3 of \cite{ParryPollicott}).
Recall that a subshift of finite type $(\Sigma_A,\s)$ is mixing if its defining matrix $A$ is aperiodic: there exists an integer $N$ such that $A^N$ has strictly positive entries.

\begin{proposition}[The Variational Principle, Theorem 3.5 \cite{ParryPollicott}]\label{prop.vp}
Let $\Sigma_A$ be a mixing subshift of finite type. Then the supremum $\sup_{\lambda \in \Mcc(\Sigma_A, \sigma)}\left\{h(\s, \lambda)\right\}$
is attained by a unique $\sigma$-invariant Borel probability measure $\m$ on $\Sigma_A$ called the measure of maximal entropy. Furthermore the above supremum agrees with the topological entropy of $(\Sigma_A, \s)$ which is obtained as the logarithm of the real, simple, maximal eigenvalue of $A$.
\end{proposition}

\subsubsection{Large deviations}
Throughout this work we will employ various statistical results to compare the length spectra of different metrics.  We will use the following lemma to establish a large deviation theorem. 

\begin{lemma}[Chapter V, Lemma 6.2 of \cite{Bougerol.Lacroix}]\label{lem.localldt}
Let $\mu_T$ be a one-parameter family of probability measures on $(X, \mathcal{B})$, where $X$ is a space and $\mathcal{B}$ is a $\sigma$-algebra on $X$. Suppose that $Z: X \to \R$ is a measurable function and that there exists $\kappa>0$ such that
\[
\lim_{T\to\infty} \frac{1}{T}\log \mathbb{E}_T(e^{tZ}) =: c(t)
\]
exists for all $t \in [-\kappa, \kappa]$. Here $\mathbb{E}_T$ represents the expectation with respect to $\mu_T$.
If $c$ is continuously differentiable on $[-\kappa,\kappa]$, then for any $\epsilon >0$
\[
\limsup_{T\to\infty} \frac{1}{T}\log \mu_T\left(x \in X : \left|\frac{Z(x)}{T} - c'(0)\right|>\epsilon\right) < 0.
\]
\end{lemma}

We present a proof since the above result is stated slightly differently than in \cite{Bougerol.Lacroix}
\begin{proof}
Let $Z_T = Z - Tc'(0)$ for $T >0$.
We note that for $\epsilon >0$ and for any $ 0 < t < \kappa$ we have that 
\[
\mu_T(Z_T> T \epsilon) \le \mathbb{E}_T(e^{t(Z_T - T\epsilon)}) = e^{-t T \epsilon - tTc'(0)} \, \mathbb{E}_T(e^{tZ})
\]
and so
\[
\limsup_{T\to\infty} \frac{1}{T} \log \mu_T(Z_T > T\epsilon) \le - t \epsilon - tc'(0) + c(t).
\]
Since $c$ is $C^1$ we can then find $t_0 \in (0,\kappa)$ such that 
\[
c(t_0) = c(t_0) - c(0) < t_0\epsilon + t_0c'(0)
\]
 and so we deduce that
\[
\limsup_{T\to\infty} \frac{1}{T} \log \mu_T(Z - Tc'(0) > T\epsilon) \le - t_0\epsilon -t_0c'(0) + c(t_0) < 0.
\]
We can apply the same argument on the interval $[-\kappa, 0)$ to obtain exponential decay when $Z -Tc'(0) < T\epsilon$. 
\end{proof}

\section{Spectrally rigid sets for Outer Space}\label{sec.outer}
In this section we prove our main results for Outer space. Recall that a hyperbolic space $(X,d)$ is called $\textnormal{CAT}(-1)$ if $(X,d)$ is geodesic and geodesic triangles in $(X,d)$ are thinner than those in the hyperbolic plane.
\begin{definition}
Fix a group $G$.
We say that a subset $E \subset \conj(G)$ is $\text{CAT}(-1)$-spectrally rigid if the following holds.  When $G$ acts cocompactly and properly discontinuously by isometries on two $\textnormal{CAT}(-1)$ metric spaces $(X_1,d_1), (X_2, d_2)$ and $\ell_{X_1}[g] = \ell_{X_2}[g]$ for all $[g] \in E$, then $\ell_{X_1}[g] = \ell_{X_2}[g]$ for all $[g] \in \conj(G)$.
\end{definition}

We begin with a crucial lemma that allows us to compare counting results for displacement functions with those for translation lengths. Throughout this section we will work with three actions of $G$ on hyperbolic metric spaces $X,X_1,X_2$. We fix base points for these spaces and abuse notation by calling them all $o$. Recall that a conjugacy class $[g] \in \conj(G)$ is called primitive if there does not exist $h \in G$ and $n > 1$ such that $h^n \in [g]$, i.e. $[g]$ is not a power.
\begin{lemma}\label{lem.hyp}
    Suppose $G$ acts properly discontinuously, cocompactly by isometries on three hyperbolic metric space $(X,d), (X_1,d_1), (X_2,d_2)$. Then there exists $C, M_1, M_2>0$ such that
    \[
    \#\{ h \in [g] : \max\{ |d(o,ho) - \ell_{X}[h]|, |d_1(o,ho) - \ell_{X_1}[h]|, |d_2(o,ho) - \ell_{X_2}[h]|\} < C  \} \ge M_1 \ell_X[g] - M_2
    \]
    for any $g \in G$ such that $[g]$ is primitive. 
\end{lemma}
\begin{proof}
Fix a word metric  $|\cdot|_S$ on $G$ associated to a finite generating set $S$. It is proved during the proof of Theorem 1.2 in Section 4.4 (page 135) of \cite{ant} that there exists $C, M_1, M_2>0$ such that
    \[
    \#\left\{ h \in [g] : |h|_S = \inf_{f \in [g]}|f|_S \right\} \ge M_1 \ell_S[g] - M_2
    \]
    for all $g \in G$ with $[g]$ primitive. Here $M_1$ is determined by
    \[
    \frac{1}{M_1} = \max\{ \#H : H \text{ is a finite subgroup of $G$}\}.
    \]
 Since there exist $C_1,C_2$ such that $||x|_S - \ell_S(x) + 2(x|x^{-1})_S| < C_1$  \cite[ Lem. 2.2 and Lem. 2.3, Ch. 9]{cdp} and $|\ell_S[x] - \inf_{y\in[x]}|y|_S| < C_2$ for all $x \in G$ \cite[Prop. 6.4, Ch. 10]{cdp}, we deduce that there is $C' >0$ such that 
     \[
    \#\{ h \in [g] : (h|h^{-1})_S < C' \} \ge M_1 \ell_S[g] - M_2
    \]
    for primitive $[g]$. Since $d,d_1,d_2$ are all pairwise quasi-isometric to $|\cdot|_S$, and this quasi-isometry inequality extends to the Gromov products (\ref{eq.gromprod}), we deduce that there exist $\widetilde{C}, \widetilde{M}_1,\widetilde{M}_2 >0 $ such that
        \[
    \#\{ h \in [g] : \max\{ (h|h^{-1})_{d}, (h|h^{-1})_{d_1}, (h|h^{-1})_{d_2}\} < \widetilde{C} \} \ge \widetilde{M}_1 \ell_X[g] - \widetilde{M}_2.
    \]
    Finally, we apply Lemma 3.2 of \cite{cantrell.tanaka.1} (which holds for the roughly geodesic pseudometrics in $\Dc_G$) to deduce the lemma.
\end{proof}

The following is a key lemma that we will apply multiple times throughout this work.

\begin{lemma}\label{lem.combo}
Suppose that $d_1, d_2, d_3 \in \Dc_G$ and for each $r \in \R$ define $d_r : G \times G \to \R$ by
\[
d_r(x,y) := rd_1(x,y) + (d_2(x,y) - d_3(x,y))  \ \ \text{ for } x,y \in G.
\]
Let $r_0$ be the critical exponent of 
\[
\sum_{x \in G} e^{- d_r(o,x)}
\]
as $r$ varies (i.e. the series converges for $r > r_0$ and diverges for $r < r_0$). Then there exist a constant $C >0 $ and a pseudometric $d \in \Dc_G$ such that
\[
d(g,h) = 
     \begin{cases}
      d_{r_0}(g,h) + C  &\text{ if $g \neq h$}\\
       0 &\text{ otherwise}
     \end{cases}
\]
belongs to $\Dc_G$. i.e. $d_{r_0}$ is within uniformly bounded distance of $d \in \Dc_G$.
\end{lemma}
\begin{proof}
Suppose that $r_0 \ge 0$. When this is the case $r_0d_1 + d_2$ is a positive linear combination of pseudometrics in $\Dc_G$ and so belongs to $\Dc_G$ by \cite[Lemma 4.1]{Oregon-Reyes}. \\

\noindent We now claim that $\theta_{(r_0d_1 + d_2)/d_3}(1) = -1.$

\textit{proof of claim: } From the definition of the Manhattan curve, it suffices to show that
\begin{equation}\label{eq.bound}
\sum_{x \in G} e^{-r_0d_1(o,x) - d_2(o,x) - sd_3(o,x)}
\end{equation}
converges when $s > -1$ and diverges when $s < -1$. Fix $s > -1$. Then 
\begin{align*}
r_0d_1(o,x) + d_2(o,x) + sd_3(o,x) &= r_0d_1(o,x) + (d_2(o,x) - d_3(o,x)) + (s+1)d_3(o,x)\\
&= d_{r_0}(o,x) + (s+1)d_3(o,x).
\end{align*}
Since $s+1 >0$ and $d_3$ is quasi-isometric to $d_1$ (as they are both in $\Dc_G$) there exist $\tau, C > 0$ such that $(s+1)d_3(o,x) > \tau d_1(o,x) -C$ for all $x \in G$.  In particular, from the previous display equation,
\begin{equation}\label{eq.metineq}
r_0d_1(o,x) + d_2(o,x) + sd_3(o,x)  = d_{r_0}(o,x) + (s+1)d_3(o,x) \ge d_{r_0 + \tau}(o,x) - C
\end{equation}
for all $x \in G$. However, by the definition of $r_0$ we know that
\[
\sum_{x \in G} e^{-d_{r_0 +\tau}(o,x)}
\]
converges. Along with $(\ref{eq.metineq})$ this implies that the series in $(\ref{eq.bound})$ converges for $s > -1$ as required. A similar argument shows that $(\ref{eq.bound})$ diverges when $s < -1$ concluding the proof of the claim.\\

Now, by this claim and  Lemma \ref{lem.bound} there is a pseudometric $d \in \Dc_G$ within uniformly bounded distance of 
\[
1(r_0d_1 + d_2)  + \theta_{(r_0d_1 + d_2)/d_3}(1) d_3 = r_0d_1 +d_2 - d_3 = d_{r_0},
\]
as required, in the case that $r_0 \ge 0$.

If $r_0 < 0$ then we can repeat the same proof but instead pair $r_0d_1$ with $d_3$. That is, $-r_0 d_1 + d_3 \in \Dc_G$, $\theta_{d_2/(-r_0d_1 + d_3)}(1) = -1$ and 
\[
1(d_2)  + \theta_{d_2/(-r_0d_1 + d_3)}(1) (-r_0d_1 + d_3) = d_{r_0} 
\]
is within uniformly bounded distance of some $d \in \Dc_G$.
\end{proof}

We also need the following large deviations result.
\begin{proposition}\label{prop.ldt}
Suppose that a group $G$ acts properly discontinuously, cocompactly by isometries on three hyperbolic metric spaces $(X, d), (X_1,d_1), (X_2,d_2)$. Then, there exists $\Lambda_X(X_1,X_2) \in \mathbb{R}$ such that for any $\epsilon > 0$, the quantity
    \[
  \frac{1}{\#\{[g] \in \conj(G): \ell_X[g] < T\}} \#\left\{ [g] \in \conj(G) : \ell_X[g] < T, \left| \frac{\ell_{X_1}[g] - \ell_{X_2}[g]}{\ell_X[g]} - \Lambda_X(X_1,X_2) \right| > \epsilon \right\}
    \]
    decays to $0$ exponentially quickly as $T\to\infty$.
\end{proposition}

\begin{proof}
Consider the series
\[
\sum_{g \in G} e^{-sd(o,go) - (d_1(o,go) - d_2(o,go))}
\]
as $s$ varies in $\R$. Let $s_0$ denote the corresponding abscissa of convergence. Then, by Lemma \ref{lem.combo} there is a pseudometric $\widetilde{d}\in \Dc_G$ that is within uniformly bounded distance of
\[
(g,h) \mapsto s_0d(go,ho) + (d_1(go,ho) - d_2(go,ho)) \, \text{ for all $g,h \in G$}.
\]
Then, by $(\ref{eq.count})$, for all sufficiently large $R>0$ and for each $t \in \R$ we have that
\[
\theta_{\widetilde{d}/d}(t) = \lim_{T\to\infty} \frac{1}{T} \log \sum_{|d(o,go) - T| <R} e^{-t\widetilde{d}(o,g)} = \lim_{T\to\infty} \frac{1}{T} \log \sum_{|\ell_X[g] - T| <R_0} e^{-t \ell_{\widetilde{d}}[g]} 
\]
where, as above, $\theta_{\widetilde{d}/d}(t)$ is the Manhattan curve of $(\widetilde{d},d)$ evaluated at $t$. In other words, if we set $ \Cc_R(T) = \{[g] \in \conj(G) : |\ell_X[g] - T| < R \}$ for $T>0$ and write
\[
\mu_T = \frac{1}{\#\Cc_T(T)} \sum_{[g] \in \Cc_R(T)} \delta_{[g]}
\]
for the uniform counting measures over $\Cc_R(T)$, then
\[
\lim_{T\to\infty} \frac{1}{T}\log \mathbb{E}_T(e^{-t\ell_{\widetilde{d}}}) = \lim_{T\to\infty} \frac{1}{T} \log \int_{\conj(G)} e^{-t\ell_{\widetilde{d}}[g]} \ d \mu_T[g] = \theta_{\widetilde{d}/d}(t) - \theta_{\widetilde{d}/d}(0).
\]
Here $\mathbb{E}_T$ is the expectation taken with respect to $\mu_T$ and we note that $\theta_{\widetilde{d}/d}(0)$ is the exponential growth rate of $\# \Cc_R(T)$ as $T\to\infty$ (see $(\ref{eq.count})$). Since $\theta_{\widetilde{d}/d}(t)$ is $C^1$ as $t$ varies in $\R$ by Theorem \ref{thm.Man},  we can then apply Lemma \ref{lem.localldt} to deduce the following. There exists $\tau \in \R$ such that for any fixed, sufficiently large $R >0$, and for any $\epsilon > 0$ the quantity
\[
  \#\left\{ [g] \in \Cc_R(T): \left| \frac{\ell_{\widetilde{d}}[g]}{T} - \tau \right| > \epsilon \right\},
  \]
grows exponentially slower than $\#\Cc_R(T)$ as $T\to\infty$. It is then a simple algebraic manipulation to verify that, for sufficiently large $R >0,$ and for any $\epsilon > 0$, the quantity
\[
  \#\left\{ [g] \in \Cc_R(T): \left| \frac{\ell_{\widetilde{d}}[g]}{\ell_X[g]} - \tau \right| > \epsilon \right\}
  \]
  grows exponentially slower than $\#\Cc_R(T)$ as $T\to\infty$. 
Now note that (since $\widetilde{d}$ is within uniformly bounded distance of the linear combination  $s_0d + d_1 - d_2$)
  \[
  \frac{\ell_{\widetilde{d}}[g]}{\ell_X[g]} = \frac{s_0\ell_X[g] + \ell_{X_1}[g] - \ell_{X_2}[g] }{\ell_X[g]} = \frac{\ell_{X_1}[g] - \ell_{X_2}[g]}{\ell_X[g]} + s_0.
  \]
From this we see that for sufficiently large $R >0$ and for any $\epsilon > 0$ the quantity
\[
  \#\left\{ [g] \in \Cc_R(T): \left| \frac{\ell_{X_1}[g] - \ell_{X_2}[g]}{\ell_X[g]} - \Lambda_X(X_1,X_2) \right| > \epsilon \right\},
  \]
  where $\Lambda_X(X_1,X_2) = \tau - s_0$ grows exponentially slower than $\Cc_R(T)$ as $T\to\infty$. It is easy to check that this implies the statement of the proposition (where the counting is over all $[g]$ with $\ell_X[g] < T$).
\end{proof}

Before we prove the main theorem of this section we prove the following simple lemma.

\begin{lemma}\label{lem.simple}
Take $T \in \cv_N$. Then the translation length function $\|\cdot\|_T$ agrees with the stable translation function, i.e.
\[
\|g\|_T = \ell_T(g) = \lim_{n\to\infty} \frac{d_T(o,g^no)}{n}
\]
for all $g \in F_N$.
\end{lemma}
\begin{proof}
Note that for any $x \in [g]$,  $\ell_T(g) = \lim_{n\to\infty} d_T(o,x^no)/n \le d_T(o,xo)$ (by the triangle inequality) and so $\ell_T(g) \le \|g\|_T$. Conversely, note that for any integer $k \ge 1$, $k\|g\|_T = \|g^k\|_T \le d_T(o,g^ko)$. Dividing this equation by $k$ and letting $k$ tend to infinity gives $\|g\|_T \le \ell_T(g)$ which proves the other inequality.
\end{proof}

We are now ready to prove our main theorem of the section.

\begin{theorem}\label{thm.cat}
Let $G$ be a finitely generated group equipped with a finite generating set $S$. If $E \subset \conj(G)$ satisfies
\[
\limsup_{T\to\infty} \sqrt{T} \, E_S(T)= \infty
\]
then $E$ is $\textnormal{CAT}(-1)$-spectrally rigid.
\end{theorem}
\begin{proof}
Suppose that $G$ acts cocompactly, properly discontinuously by isometries of two $\textnormal{CAT}(-1)$ metric spaces $(X_1,d_1), (X_2,d_2)$ and that $S$ is a fixed finite generating set on $G$. Suppose that $d_1$ and $d_2$ have different marked length spectra, we will show that 
\[
\limsup_{T\to\infty} \frac{\sqrt{T} \, \, \#\{ [g] \in \conj(G): \ell_{X}[g] < T, \ell_{X_1}[g] = \ell_{X_2}[g]\}}{\#\{ [g] \in \conj(G): \ell_{\mathfrak{g}}[g] < T \}} < \infty.
\]
We consider the asymptotic average
\[
\Lambda_S(X_1,X_2) :=  \lim_{T\to\infty} \frac{1}{\#\{[g] \in \conj(G): \ell_S[g] < T\}} \sum_{\ell_S[g] < T} \frac{\ell_{X_1}[g] - \ell_{X_2}[g]}{\ell_S[g]}
\]
which exists by Proposition \ref{prop.ldt}. We now consider two cases depending on the value of $\Lambda_S(X_1,X_2)$.\\

    \textit{Case 1: When $\Lambda_S(X_1,X_2) \neq 0$.} Then by Proposition \ref{prop.ldt} for any $\epsilon > 0$, the count
    \[
   \#\left\{ [g] \in \conj(G) : \ell_S[g] < T, \left| \frac{\ell_{X_1}[g] - \ell_{X_2}[g]}{\ell_S[g]} - \Lambda_S(X_1,X_2) \right| > \epsilon \right\} 
    \]
    grows strictly exponentially slower than $\#\{[g] \in \conj(G): \ell_S[g] < T\}$. By the assumption that $ \Lambda_S(X_1,X_2)  \neq 0$, we deduce that 
    \[
    \frac{\#\{[g] \in \conj(F_N) : \ell_S[g] < T, \ell_{X_1}[g] = \ell_{X_2}[g]\}}{\#\{[g] \in \conj(F_N) : \ell_S[g] < T\}}
    \]
    decays exponentially quickly to $0$ as $T\to\infty$. \\
    
    \textit{Case 2: When $\Lambda_S(X_1,X_2) = 0$.}
    Let $\Cc_S(T) = \{ [g] \in \conj(G) : \ell_S[g] < T\}$ for $T > 0$ and observe that
    \begin{equation}\label{eq.conj}
    \frac{1}{\#\mathcal{C}_S(T)} \#\left\{ [g] \in \mathcal{C}_S(T) : \ell_{X_1}[g] = \ell_{X_2}[g] \right\} 
    \end{equation}
      is bounded above by
\[
     \frac{1}{\#\mathcal{C}_S(T)} \#\left\{ [g] \in \mathcal{C}_S(T) : \frac{\ell_{X_1}[g] - \ell_{X_2}[g]}{\sqrt{T}} \in \left(-\frac{1}{\sqrt{T}}, \frac{1}{\sqrt{T}}\right) \right\}
\]
for each $T >0$.
Since the number of non primitive conjugacy classes of $\ell_S$ length less than $T$ grows exponentially slower than $\mathcal{C}_S(T)$, it suffices to show that 
    \begin{equation}\label{eq.prim}
    \frac{1}{\#\mathcal{C}_S(T)}\#\left\{ [g] \in \mathcal{C}_S(T) : \textnormal{$[g]$ is primitive, $\ell_S[g] > T/2$ and } \, \frac{\ell_{X_1}[g] - \ell_{X_2}[g]}{\sqrt{T}} \in \left(-\frac{1}{\sqrt{T}}, \frac{1}{\sqrt{T}}\right) \right\} 
    \end{equation}
    is $O\left( T^{-1/2}\right)$    as $T\to\infty$. To do this we leverage statistical results for the displacement functions $d_1, d_2$. First we use Lemma \ref{lem.hyp} to deduce that if $[g]$ is primitive and satisfies $\ell_S(g) > T/2$ then there are at least $M_1T/2 - M_2$ elements $x \in [g]$ satisfying
     \[
     |g|_S < T + C, \ \  |\ell_{X_1}[g] - d_1(o,go)| < C \ \ \text{and} \ \   |\ell_{X_2}[g] - d_2(o,go)| < C.
    \]
      Hence there exist $C_1,C_2, C_3 > 0$ such that the right hand side of $(\ref{eq.prim})$ is bounded above by
    \[
    C_1\left(\frac{M_1 T}{2} - M_2\right)^{-1} \#\left\{ g \in S(T + C_2) : \frac{d_1(o,go) - d_2(o,go)}{\sqrt{T+C_2}} \in \left( -\frac{C_3}{\sqrt{T}}, \frac{C_3}{\sqrt{T}} \right)   \right\}
    \]
    for all $T >0$ sufficiently large where $S(T) = \{ x \in G : |x|_S < T\}$. We now use the fact that there exists $C>0$ such that
    \[
    \frac{1}{T\#\mathcal{C}_S(T)} \le \frac{C}{\# S(T)}
    \]
  for $T >0$ sufficiently large  \cite{CoornaertKnieper} to see that there is $\widetilde{C}_1,\widetilde{C}_2 >0$ such that $(\ref{eq.conj})$ is bounded above by
    \[
 \, \frac{ \widetilde{C}_1}{\#S(T 
    + \widetilde{C}_2)} \#\left\{ g \in S(T + \widetilde{C}_2) : \frac{d_1(o,go) - d_2(o,go)}{\sqrt{T + \widetilde{C}_2}} \in \left( -\frac{\widetilde{C}_2}{\sqrt{T}} , \frac{\widetilde{C}_2}{\sqrt{T}} \right) \right\}.
    \]
    By the central limit theorem with Berry-Esseen error term \cite[Theorem 1.2]{stats} (which is stated for counting over spheres but also holds for counting over balls with the same proof) for the displacement functions we have that
    \[
    \frac{1}{\#S(T 
    + \widetilde{C})} \#\left\{ g \in S(T + \widetilde{C}) : \frac{d_1(o,go) - d_2(o,go)}{\sqrt{T + \widetilde{C}}} \in \left( -\frac{\widetilde{C}}{\sqrt{T}} , \frac{\widetilde{C}}{\sqrt{T}} \right) \right\} = O\left( \frac{1}{\sqrt{T}}\right)
    \]
    as $T \to \infty$.
    
    This concludes the proof of the theorem. 
\end{proof}

\begin{remark}
The line of argument that we used in the above theorem actually shows a stronger result: if $E \subset \conj(G)$ for which there exist $C >0$ and $\theta \in [0,1/2)$ such that 
\[
|\ell_{X_1}[g] - \ell_{X_2}[g]| \le C \ell_S[g]^\theta \ \ \text{ for all $[g]\in E$}
\]
then $E_S(T) = O(T^{\theta - 1/2})$ as $T\to\infty$.
\end{remark}

We can now prove our main theorem for Outer Space.

\begin{proof}[Proof of Theorem \ref{thm.outermain}]
The first part of the theorem now follows immediately from Theorem \ref{thm.cat} and Lemma \ref{lem.simple}. We establish the optimality of the theorem in the proposition below.
\end{proof}

\begin{proposition}\label{prop.optimal}
Take $N \ge 3$ and let $A$ denote the canonical free generating set for $F_N$. Then there are $T_1, T_2 \in \cv_N$ for which there exists $C>0$ with
    \[
    \frac{ \# \{ [g] \in \textnormal{conj}(F_N): \ell_A[g] = n, \|g\|_{T_1}= \|g\|_{T_2}\}}{\#\{[g] \in \textnormal{conj}(F_N): \ell_A[g] = n\}} \sim \frac{C}{\sqrt{n}}
    \]
    as $n\to\infty$. When $N=2$ there are $T_1, T_2 \in \cv_2$ and $C'>0$ such that, for the canonical free generating set $A$, 
    \[
      \frac{ \# \{ [g] \in \textnormal{conj}(F_N): \ell_A[g] = 2n, \|g\|_{T_1}= \|g\|_{T_2}\}}{\#\{[g] \in \textnormal{conj}(F_N): \ell_A[g] = 2n\}} \sim \frac{C'}{\sqrt{n}}
    \]
    as $n\to\infty$.
\end{proposition}

Here the notation `$\sim$' means the following: if $a_n, b_n$ are sequences we write $a_n \sim b_n$ if $a_n/b_n \to 1$ as $n\to\infty$.

\begin{proof}
Consider the bouquet of $N \ge 2$ loops. A labelling for the bouquet of $N$ circles is a choice of length for each loop.  
Consider the following  labellings.  Labelling $1$ assigns each of the $N$ loops the  length $2$. Labelling $2$ assigns $N-2$ of the loops length $2$ and the other $2$ loops are assigned lengths $1$ and $3$.
These two labellings give rise to (cocompact, properly discontinuous) actions of the free group $F_N$ on an $\mathbb{R}$ tree (the universal covers of the labelled bouquets).
These actions correspond to points $T_1, T_2 \in \text{cv}_N$ in Outer Space.
Our choice of labellings guarantees that $T_1 \neq T_2$. Indeed, by construction we can find generators $a,b \in A$ such that $\|a\|_{T_1} = 2, \|a\|_{T_2} = 1$ and $\|b\|_{T_2} = 2, \|b\|_{T_2}=3$.  This guarantees $T_1, T_2$ determine different points in $\text{cv}_N$. 

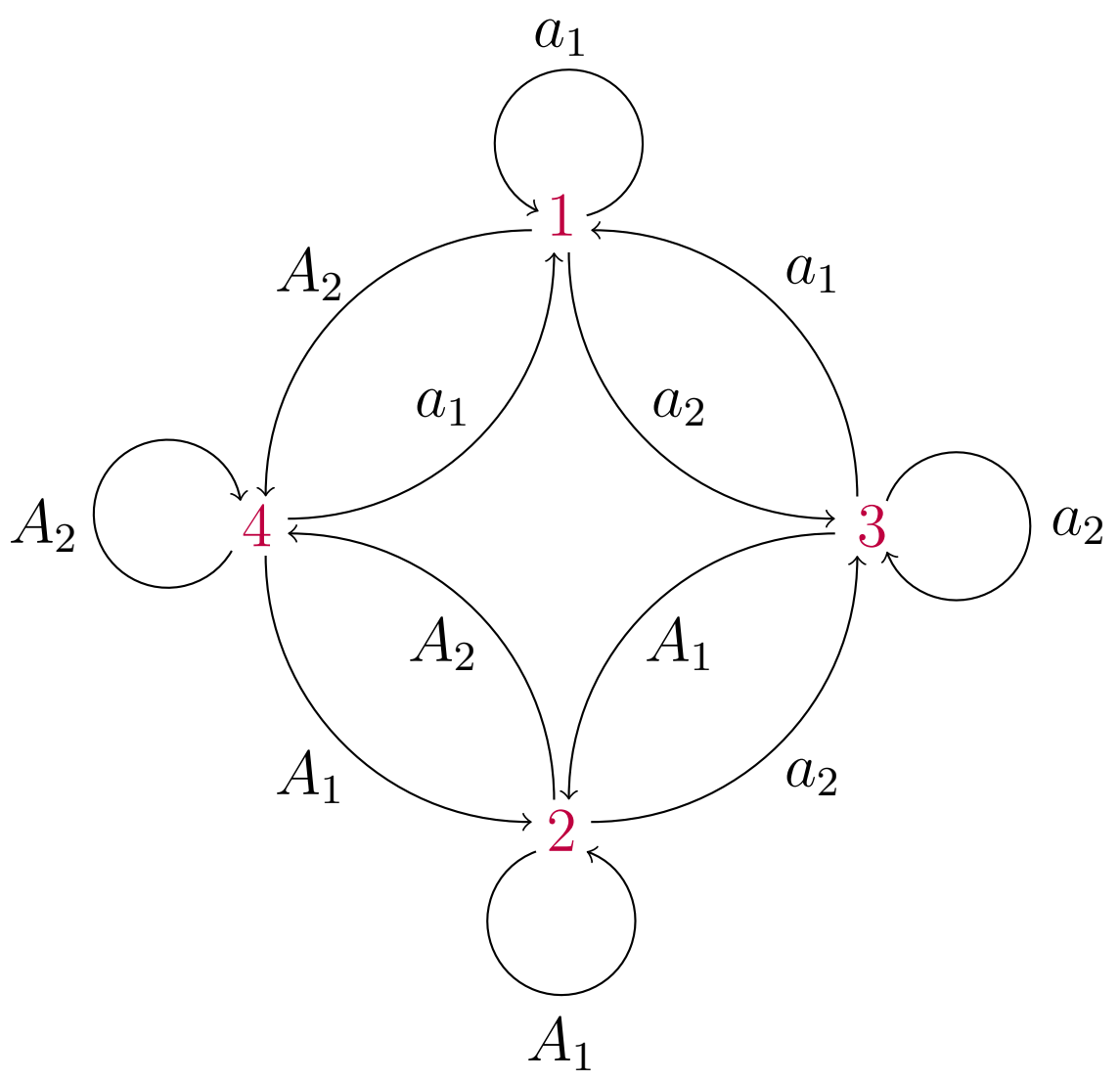
\begin{figure}[h!tbp]
\begin{tikzpicture}
\draw[purple] (0,2.1)  node {$1$};
\draw[purple] (2.1,0) node {$3$};
\draw[purple] (0,-2.06) node {$2$};
\draw[purple] (-2.06,0)  node {$4$};
\draw[->] (2,0.2) arc (0:90:1.8cm);
\draw[->] (-0.05,-1.85) arc (0:90:1.8cm);
\draw[<-] (1.85,0.05) arc (270:180:1.8cm);
\draw[<-] (-0.2,-2) arc (270:180:1.8cm);
\draw[->] (-0.2,2) arc (90:180:1.8cm);
\draw[->] (1.85,-0.05) arc (90:180:1.8cm);
\draw[->] (0.2,-2) arc (-90:0:1.8cm);
\draw[->] (-1.85,0.05) arc (-90:0:1.8cm);
\draw[->] (0.17,2.1) arc (-76:246:0.5cm);
\draw[->] (2.2,0.17) arc (160: - 160:0.5cm);
\draw[<-] (0.17,-2.2) arc ( 70: -250: 0.5cm);
\draw[<-] (-2.17, 0.17) arc ( 10: 330: 0.5cm);
\draw (3.5,0)  node {$a_2$};
\draw (-3.5,0)  node {$A_2$};
\draw (0, 3.3)  node {$a_1$};
\draw (0,-3.5)  node {$A_1$};
\draw (-1.7,-1.7)  node {$A_1$};
\draw (-1.7,1.7)  node {$A_2$};
\draw (1.7,-1.7)  node {$a_2$};
\draw (1.7,1.7)  node {$a_1$};
\draw (0.8,0.8)  node {$a_2$};
\draw (-0.8,0.8)  node {$a_1$};
\draw (0.8,-0.8)  node {$A_1$};
\draw (-0.8,-0.8)  node {$A_2$};
\end{tikzpicture}
\caption{The graph described by $A_2$.}
\end{figure}

Now consider the subshift of finite type $\Sigma_N$ defined using the $2N \times 2N$ transition matrix $A_N$ with entries
\[
A_N(i,j) =   \begin{cases}
     0 &\text{if $i$ is odd and $j=i+1$, or if $i$ is even and $j=i-1$ }\\
             1 &\text{otherwise.}
     \end{cases}
\]
Consider the finite directed graph induced by $A$: let $\mathcal{G}$ be the graph with $2N$ vertices labelled $1,\ldots, 2N$ such that vertices $i,j$ are connected with a directed edge from $i$ to $j$ if and only if $A(i,j) =1$. Figure 1 depicts the directed graph determined by the matrix $A_2$.

Let $\{a_1,a_1^{-1},\ldots, a_N,a_N^{-1}\}$ be the canonical free generating set of $F_N$ and label the edges of $\mathcal{G}$ in the following way: the directed edge from $i$ to $j$ is labelled by $a_{(j+1)/2}$ if $j$ is odd and $a_{j/2}^{-1}$ if $j$ is even. Figure 1 has this labelling where $A_1= a_1^{-1}$ and $A_2 = a_2^{-1}$. One can check that this labelling induces a bijection between the periodic orbits in $\Sigma_N$ and the conjugacy classes of $F_N$. Given a periodic orbit $\textnormal{Orb}(x) = \{x,\sigma(x), \ldots, \sigma^{n-1}(x)\}$ (with $\sigma^n(x) = x$) we write $[g_x]$ for the corresponding conjugacy class.

There are simple functions $r_1, r_2$ on $\Sigma_A$ corresponding to $T_1, T_2$. Indeed, define $r_1 : \Sigma_N \to \mathbb{Z}$ to be the function that is constant and equal to $2$ and define $r_2: \Sigma_N \to \Z$ by
\[
r_2(x) = \begin{cases}
     2 &\text{if $x_1 = 1,2,\ldots, 2N-4$}\\
     1 &\text{if $x_1 = 2N-3, 2N-2$}\\
     3 &\text{if $x_1 = 2N-1, 2N$}
\end{cases}
\] 
when $x = (x_n)_{n=0}^\infty$.
By construction these potentials are constant on cylinders of length $2$.  Given a periodic point $x\in\Sigma_N$ with $\sigma^n(x) = x$ we have that (by construction)
\begin{equation}\label{eq.birk}
r_1^n(x) = \ell_{T_1}[g_x] \ \ \text{ and } r_2^n(x) = \ell_{T_2}[g_x]
\end{equation}
for any periodic point $x \in \Sigma_N$. See Figure $2$ for a visual representation of the functions/labellings $r_1,r_2$ when $N=2$.

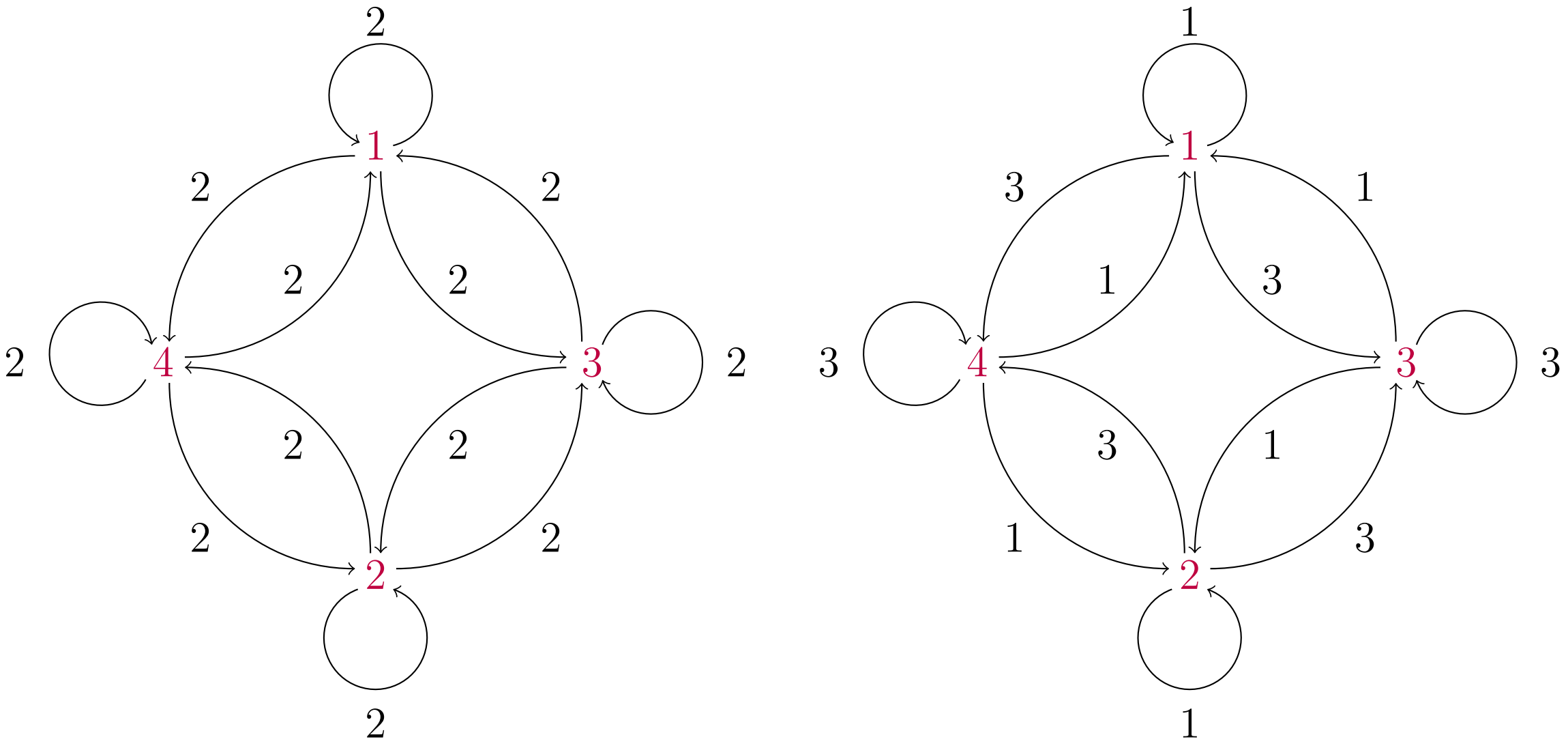
\begin{figure}[h!tbp]
\begin{tikzpicture}
\draw[purple] (0,2.1)  node {$1$};
\draw[purple] (2.1,0) node {$3$};
\draw[purple] (0,-2.06) node {$2$};
\draw[purple] (-2.06,0)  node {$4$};
\draw[->] (2,0.2) arc (0:90:1.8cm);
\draw[->] (-0.05,-1.85) arc (0:90:1.8cm);
\draw[<-] (1.85,0.05) arc (270:180:1.8cm);
\draw[<-] (-0.2,-2) arc (270:180:1.8cm);
\draw[->] (-0.2,2) arc (90:180:1.8cm);
\draw[->] (1.85,-0.05) arc (90:180:1.8cm);
\draw[->] (0.2,-2) arc (-90:0:1.8cm);
\draw[->] (-1.85,0.05) arc (-90:0:1.8cm);
\draw[->] (0.17,2.1) arc (-76:246:0.5cm);
\draw[->] (2.2,0.17) arc (160: - 160:0.5cm);
\draw[<-] (0.17,-2.2) arc ( 70: -250: 0.5cm);
\draw[<-] (-2.17, 0.17) arc ( 10: 330: 0.5cm);
\draw (3.5,0)  node {$2$};
\draw (-3.5,0)  node {$2$};
\draw (0, 3.3)  node {$2$};
\draw (0,-3.5)  node {$2$};
\draw (-1.7,-1.7)  node {$2$};
\draw (-1.7,1.7)  node {$2$};
\draw (1.7,-1.7)  node {$2$};
\draw (1.7,1.7)  node {$2$};
\draw (0.8,0.8)  node {$2$};
\draw (-0.8,0.8)  node {$2$};
\draw (0.8,-0.8)  node {$2$};
\draw (-0.8,-0.8)  node {$2$};
\end{tikzpicture}
\ \ 
\begin{tikzpicture}
\draw[purple] (0,2.1)  node {$1$};
\draw[purple] (2.1,0) node {$3$};
\draw[purple] (0,-2.06) node {$2$};
\draw[purple] (-2.06,0)  node {$4$};
\draw[->] (2,0.2) arc (0:90:1.8cm);
\draw[->] (-0.05,-1.85) arc (0:90:1.8cm);
\draw[<-] (1.85,0.05) arc (270:180:1.8cm);
\draw[<-] (-0.2,-2) arc (270:180:1.8cm);
\draw[->] (-0.2,2) arc (90:180:1.8cm);
\draw[->] (1.85,-0.05) arc (90:180:1.8cm);
\draw[->] (0.2,-2) arc (-90:0:1.8cm);
\draw[->] (-1.85,0.05) arc (-90:0:1.8cm);
\draw[->] (0.17,2.1) arc (-76:246:0.5cm);
\draw[->] (2.2,0.17) arc (160: - 160:0.5cm);
\draw[<-] (0.17,-2.2) arc ( 70: -250: 0.5cm);
\draw[<-] (-2.17, 0.17) arc ( 10: 330: 0.5cm);
\draw (3.5,0)  node {$3$};
\draw (-3.5,0)  node {$3$};
\draw (0, 3.3)  node {$1$};
\draw (0,-3.5)  node {$1$};
\draw (-1.7,-1.7)  node {$1$};
\draw (-1.7,1.7)  node {$3$};
\draw (1.7,-1.7)  node {$3$};
\draw (1.7,1.7)  node {$1$};
\draw (0.8,0.8)  node {$3$};
\draw (-0.8,0.8)  node {$1$};
\draw (0.8,-0.8)  node {$1$};
\draw (-0.8,-0.8)  node {$3$};
\end{tikzpicture}
\caption{The labellings corresponding to the functions $r_1$ (left) and  $r_2$ (right) when $N =2$.}
\end{figure}

Let $\mu$ denote the measure of maximal entropy of $\Sigma_N$ (see Proposition \ref{prop.vp}). Since the $r_j$ are constant on $2$ cylinders for $j=1,2$ we have that
\begin{equation}\label{eq.int}
\int_{\Sigma_N} r_j \ d\mu = \sum_{x_0, x_1 :  A(x_0,x_1) = 1} l_j(x_0,x_1) \mu[x_0,x_1]    
\end{equation}
where
\begin{itemize}
\item the sum is over all $x_0, x_1 \in \{1,\ldots, 2N\}$ with $A(x_0,x_1) = 1$;
\item  $l_j(x_0,x_1)$ is the length of the edge connecting $x_0$ to $x_1$ with respect to Labelling $j$ ($j=1,2$); and,
\item $\mu[x_0,x_1]$ is the measure of the cylinder set $[x_0,x_1]$.
\end{itemize}
 We can explicitly compute the measure of maximal entropy. By a result of Parry \cite{parry}, the measure of maximal entropy is 
 defined as follows: the $\mu$ measure of any cylinder set $[x_0,x_1,\ldots,x_{n-1}] \subset \Sigma_N$ is $\frac{1}{2N(2N-1)^{n-1}}$ 
 (and this uniquely extends to a $\sigma$-invariant measure on $\Sigma_N$).
 Hence for any $x_0,x_1$ with $A(x_0,x_1) = 1$, $\mu[x_0,x_1] = \frac{1}{2N(2N-1)}$. It follows from this and $(\ref{eq.int})$ that
\[
\int_{\Sigma_N} r_1 \ d\mu = \int_{\Sigma_N} r_2 \ d\mu = \frac{2\#\{ x_0,x_1 \in \{1,\ldots, 2N\}: A(x_0,x_1) =1\}}{2N(2N-1)} = 2.
\]
Now consider the function $R: \Sigma_N \to \mathbb{R}$ defined by $R = r_1 - r_2$. This function is constant on $2$ cylinders and takes values in $\mathbb{Z}$. By construction and the above discussion we know that
\begin{enumerate}
\item $R$ has average $0$ with respect to the measure of maximal entropy:
\[
\int_{\Sigma_N}R \ d\mu = \int_{\Sigma_N} r_1 \ d\mu - \int_{\Sigma_N} r_2 \ d\mu =  0
\]
where $\mu$ is the measure of maximal entropy; and,
\item for each $m \in \mathbb{Z}$ there exists a periodic point $x \in \Sigma_N$ with $\sigma^k(x) = x$ such that $R^k(x) = m$.
\end{enumerate}
These conditions allow us to apply a result of Pollicott and Sharp \cite[Theorem 1]{pollicott.sharp.rates} to deduce that there exists a constant $C >0$ and integer $k \in \Z_{\ge 0}$ such that
\begin{equation}\label{eq.k}
\#\{x \in \Sigma_N : \sigma^{kn}(x) = x \text{ and } R^{kn}(x) = 0 \} \sim \frac{C \, (2N-1)^{kn}}{\sqrt{n}}
\end{equation}
as $n\to\infty$. Here we have used the fact that the topological entropy of $(\Sigma_A,\sigma)$ is $h = \log(2N-1)$: it is easy to see that the real, simple maximal eigenvalue of $A_N$ is $2N-1$. Furthermore the integer $k$ is the cardinality of the quotient $|\Z/\Delta_R|$ where $\Delta_R$ is the subgroup generated by the values
\[
\bigcup_{n\ge 1} \{ R^n(x) - R^n(y) : x,y \in \Sigma_N, \sigma^n(x) = x, \sigma^n(y) = n\}.
\]
It is easy to check that $k=1$ when $N \ge 3$ and $k=2$ when $N=2$. This leads to the two cases in the theorem.

Suppose that $N\ge 3$ so that $k=1$ in $(\ref{eq.k})$.
We claim that  $(\ref{eq.k})$ implies that
\begin{equation}\label{eq.asym}
\#\{ [g] \in \text{conj}(F_N) : \ell_A[g] = n, \, \ell_1[g] = \ell_2[g] \} \sim \frac{C \, (2N-1)^n}{n^{3/2}}
\end{equation}
for some $C >0$. To see this note the following facts:
\begin{enumerate}
    \item $R^n(x) = \ell_1[g_x] - \ell_2[g_x]$ by $(\ref{eq.birk})$;
    \item for every periodic orbit of least period $n$ there are $n$ corresponding periodic points; and,
    \item the number of periodic points of non-least period $n$ grows exponentially slower than the number of periodic points of least period $n$, i.e.
    \[
    \#\{ x \in \Sigma_N: \sigma^n(x) = x \text{ and $\sigma^k(x) = x$ for some $k < n$} \} = o((2N-1 - \delta)^n)
    \]
    for some $\delta >0$  as $n\to\infty$.
\end{enumerate}
The proof is then complete for $N\ge 3$ after using the fact that
\[
\#\{ [g] \in \text{conj}(F_N) : \ell_A[g] = n\} \sim \frac{(2N-1)^n}{n}
\]
as $n\to\infty$ and dividing equation $(\ref{eq.asym})$ by this quantity.

The case $N=2$ follows similarly from $(\ref{eq.k})$ with $k=2$ and the above observations.
\end{proof}

Before moving on to our study of manifolds we note the following corollary.

\begin{corollary}
Take $\phi \in \textnormal{Aut}(F_N)$, $g \in F_N$ and let $\Sigma = \{ [\phi^n(g)]: n \in \mathbb{Z}\} \subset \conj(F_N)$. Then for any finite generating set $S$ of $F_N$, 
\[
\Sigma_S(T) = O\left( \frac{1}{\sqrt{T}}\right) 
\]
as $T\to\infty$.
\end{corollary}
\begin{proof}
By \cite{ray.2} the set $\Sigma$ is not spectrally rigid and so by Theorem \ref{thm.outermain} its density must decay at least like $1/\sqrt{T}$.
\end{proof}


\section{Spectrally rigid sets for negatively curved manifolds}\label{sec.manifolds}

In this section we study spectrally rigid sets for closed, negatively curved manifolds. 
We  begin by noting the following corollary of Theorem \ref{thm.cat}.
\begin{corollary}\label{cor.wordlsr}
Suppose that $M$ is a closed manifold that admits smooth, negatively curved Riemannian metrics. Equip $G = \pi_1(M)$ with a finite generating set $S$. Then,  if $E \subset \conj(G)$ is a subset of closed geodesics on $M$ such that 
\[
\limsup_{T\to\infty} \frac{\sqrt{T} \, \, \#\{ [g] \in E: \ell_{S}[g] < T\}}{\#\{ [g] \in \conj(G): \ell_{S}[g] < T\}} = \infty,
\]
then $E$ is spectrally rigid.
\end{corollary}
We now turn our attention to proving Theorem \ref{thm.manifolds}. Before doing so, we need various preliminary results. Let $\gfr,\gfr_1,\gfr_2$ be negatively curved, smooth Riemannian metrics on the closed manifold $M$. 
We begin by defining 
\[
 \Lambda_\gfr(\gfr_1,\gfr_2) = \lim_{T\to\infty} \frac{1}{\#\{[g] \in \conj(G) : \ell_\gfr[g] < T\}} \sum_{[g] : \ell_\gfr[g] < T} \frac{\ell_{\gfr_1}[g] - \ell_{\gfr_2}[g]}{\ell_\gfr[g]},
\]
which exists by Proposition \ref{prop.ldt}. 

Throughout this section our results will assume different combinations of the following conditions:
\begin{enumerate}
\item $\Lambda_\gfr(\gfr_1,\gfr_2) = 0$; 
\item  we have that
\[
\inf_{[g] \in \conj(G)'} \frac{\ell_{\gfr_2}[g]}{\ell_{\gfr_1}[g]} < 1 <  \sup_{[g] \in \conj(G)'} \frac{\ell_{\gfr_2}[g]}{\ell_{\gfr_1}[g]} 
\]
where $\conj(G)'$ are the non-torsion conjugacy classes; and,
\item  $\gfr_1$ and $\gfr_2$ have distinct length spectra.
\end{enumerate}

\begin{lemma}\label{lem.nonarith}
Let $\gfr,\gfr_1,\gfr_2$ be negatively curved, smooth Riemannian metrics on the closed manifold $M$. Suppose that $(1)$ and $(3)$ hold. Then, there do not exist $a,b \in \R$, other than $a = b= 0$, such that
\[
a (\ell_{\gfr_1}[g] - \ell_{\gfr_2}[g]) + b \ell_\gfr[g] \in \mathbb{Z}
\]
for all $[g] \in \conj(G)$.
\end{lemma}
\begin{proof}
 Suppose that there exist $a,b \in \R$, not both identically $0$ such that  $a (\ell_{\gfr_1}[g] - \ell_{\gfr_2}[g]) + b \ell_\gfr[g] \in \mathbb{Z}$  for all $[g] \in \conj(G)$. Now lift $\gfr_1, \gfr_2, \gfr$ to $G$ to obtain pseudometrics that are quasi-isometric to word metrics and belong to $\Dc_G$. We call these pseudometrics $d_\gfr, d_{\gfr_1}, d_{\gfr_2}$ respectively. Let $\beta_\gfr : G \times \partial G \to \R$ be the Busemannn function associated to $d_\gfr$. The corresponding Busemannn functions for $d_{\gfr_1}$ and $d_{\gfr_2}$ are defined in the same way. As discussed in the preliminaries, these Busemannn functions are H\"older continuous cocycles: for each fixed $g \in G$, $\beta_\gfr(g,\xi)$ is H\"older in the variable $\xi$ for the natural metric on $\partial G$ and $\beta_\gfr$ satisfies the cocycle identity.
Fix a hyperbolic element $g\in G$.
The cocycle $c(g,\xi) := a(\beta_{\gfr_1}(g,\xi) - \beta_{\gfr_2}(g,\xi) )  + b\beta_\gfr(g,\xi) $ is H\"older continuous and $\partial G$ is connected. Hence we can apply Lemma 3.3 from \cite{GMM2018} to find a geodesic connecting $g^+$ and $g^-$ which is invariant under $g^n$ for some $n \ge 1$. Applying the Livsic Theorem argument in Proposition 3.2 from \cite{GMM2018} we deduce that $c(g^n,g^+) = c(g^n, g^-)$ which implies that 
\[
a (\ell_{\gfr_1}[g^n] - \ell_{\gfr_2}[g^n]) + b \ell_\gfr[g^n]  = c(g^n,g^+) = c(g^n,g^-)=  - a (\ell_{\gfr_1}[g^{-n}] - \ell_{\gfr_2}[g^{-n}]) - b \ell_\gfr[g^{-n}],
\]
where the first and last equalities follow from the cocycle identity.
It follows that $a (\ell_{\gfr_1}[g] - \ell_{\gfr_2}[g]) + b \ell_\gfr[g] = 0$ and this holds for all $[g]$. Note that if $a = 0$ then we necessarily have $b =0$ (since the length spectrum of $\gfr$ is not trivial). Assume $a \neq 0$. When this is the case $b = a\Lambda(\gfr_1, \gfr_2) = 0$ by the definition of $\Lambda(\gfr_1,\gfr_2)$. Hence we've deduced that $\ell_{\gfr_1}[g] - \ell_{\gfr_2}[g] = 0$ for all $g \in G$. This contradicts the assumption that $\gfr_1,\gfr_2$ have distinct marked length spectra and so concludes the proof.
\end{proof}

We also need the following lemma. 
\begin{lemma}\label{lem.pro}
Suppose that $(1)$ holds but that $(2)$  does not hold.
Then there exists $a,b \in \R$, not both zero, such that
\[
a\ell_\gfr[g] + b(\ell_{\gfr_1}[g] - \ell_{\gfr_2}[g]) = 0
\]
for all $[g] \in \conj(G)$.
\end{lemma}
\begin{proof}
Since $(2)$ does not hold, either $\ell_{\gfr_1}[g] - \ell_{\gfr_2}[g] \ge 0$ for all $[g] \in \conj(G)$  or $\ell_{\gfr_2}[g] - \ell_{\gfr_1}[g] \ge 0$ for all $[g] \in \conj(G)$. Without loss of generality, suppose that $\ell_{\gfr_1}[g] - \ell_{\gfr_2}[g] \ge 0$ for all $[g] \in \conj(G)$. 

Let $s_0 \in \R$ be the abscissa of convergence of the sum
\[
\sum_{g \in G} e^{-sd_\gfr(o,g) - (d_{\gfr_1}(o,g) - d_{\gfr_2}(o,g))}
\]
as $s$ varies. Here $d_\gfr, d_{\gfr_1}, d_{\gfr_2}$ are the pseudometrics introduced in the previous lemma. Then, by Lemma \ref{lem.combo} there exists $C >0$ such that
\[
d(g,h) = 
     \begin{cases}
       s_0d_\gfr(g,h) + (d_{\gfr_1}(g,h) - d_{\gfr_2}(g,h)) + C  &\text{ if $g \neq h$}\\
       0 &\text{ otherwise}
     \end{cases}
\]
belongs to $\Dc_G$. Also, by Lemma 4.1 in \cite{Oregon-Reyes} for any $\alpha > 0$ the sum $d + \alpha d_\gfr$ is a pseudometric belonging to $\Dc_G$. By adding a positive multiple of $d_\gfr$ if necessary, assume that $s_0 >0$ in the definition of $d$. 

Now consider the pseudometrics $d$ and $d_\ast: = s_0 d_\gfr$ in $\Dc_G$ and note that $\ell_d[g] \ge \ell_{d_\ast}[g]$ for all $[g] \in \conj(G)$. It follows that $v_d \le v_{d_\ast}$. Now, by Theorem \ref{thm.Man} we have that
\[
\tau(d/d_\ast):= \lim_{T\to\infty} \frac{1}{\#\{[g] \in \conj(G): \ell_{d_\ast}[g] < T\}} \sum_{[g]: \ell_{d_\ast}[g] < T} \frac{\ell_{d}[g]}{\ell_{d_\ast}[g]} \ge \frac{v_{d_\ast}}{v_d} \ge 1
\]
and $\tau(d/d_\ast)= v_{d_\ast}/v_d$ if and only if  $\ell_{d}$ is proportional to $\ell_{d_\ast}$.

Note however that assumption $(1)$ implies that $\tau(d/d_\ast) = 1$. This is because
\[
\frac{\ell_d[g]}{\ell_{d_\ast}[g]} = 1 + \frac{\ell_{\gfr_1}[g] - \ell_{\gfr_2}[g]}{s_0\ell_\gfr[g]} \ \ \text{ for all $[g] \in \conj(G)$.}
\]
From the above inequality it follows that $\tau(d/d_\ast) = v_{d_\ast}/v_d$ and so $\ell_{d}$ is proportional to $\ell_{d_\ast}$. This implies the conclusion of the lemma.
\end{proof}

Combining the previous two lemmas we deduce the following corollary.
\begin{corollary}\label{cor.imp}
Suppose that conditions $(1)$ and $(3)$ hold, then condition $(2)$ is also satisfied.
\end{corollary}

We need the following proposition.

\begin{proposition}\label{prop.asymp}
Suppose that $d,d_\ast \in \Dc_G$ have H\"older continuous Busemann cocycles that are obtained as limits as in equation $(\ref{eq.limit})$ and that condition $(2)$ holds.
Also assume that the only $a,b \in \R$ such that $a \ell_d + b \ell_{d_\ast}$ takes values in $\Z$ are $a= b= 0$.
Then there exists $\alpha > 0$ such that for any $\epsilon >0$ we can find $C>0$ with
\[
\#\{ g \in G: d_\ast(o,g) < T, |d(o,g) - d_\ast(o,g)| < \epsilon \} \sim \frac{Ce^{\alpha T}}{\sqrt{T}}
\]
as $T\to\infty$.
\end{proposition}
\begin{proof}
This is similar to the statement of \cite[Theorem 1.10]{cantrell.new} except for the following changes. In \cite[Theorem 1.10]{cantrell.new}  it is assumed that  
\begin{enumerate}[(i)]
\item $d,d_\ast$ are strongly hyperbolic; and,
\item $d,d_\ast$ have exponential growth rate $1$.
\end{enumerate} 
In the current theorem condition (i) is replaced with the assumption that $d,d_\ast$ have H\"older continuous Busemannn cocycles and (ii) is replaced by assumption $(2)$.

We note that when (i) holds, $d,d_\ast$ have H\"older continuous Busemannn cocycles that are obtained as limits as in $(\ref{eq.limit})$. Further by inspecting the proof \cite[Theorem 1.10]{cantrell.new}, strong hyperbolicity is used only to guarantee these properties. Hence we can replace  (i) with the assumption that $d,d_\ast$ have H\"older continuous Busemannn cocycles that are obtained as in $(\ref{eq.limit})$.Busemann

We will now show that it is possible to swap condition (ii) with assumption (2). To do so,  consider the pseudometric $d_\ast' = v_{d_\ast}d_\ast$ (i.e. scale $d_\ast$ so it has exponential growth rate $1$). Then, consider the Manhattan curve  $\theta_{d/d_\ast'}(t)$ for $d,d_\ast'$ (see Proposition \ref{prop.ldt}). By Lemma \ref{lem.bound} for each $t \in \R$ there exists $C_t >0$ such that we can find a pseudometric $d_t \in \Dc_G$ that is equal to 
\[
d_t(g,h) = 
     \begin{cases}
       td(g,h) + \theta_{d/d_\ast'}(t)d_\ast'(g,h)+ C_t  &\text{ if $g \neq h$}\\
       0 &\text{ otherwise}
     \end{cases}
\]
and has exponential growth rate $1$. Note that $d_t, d_\ast '$ have H\"older continuous Busemann cocycles (and have exponential growth rates $1$) that are obtained as limits as in $(\ref{eq.limit})$ and  so, as we discussed above, we can apply \cite[Theorem  1.20]{cantrell.new} to the pair $d_\ast' + C_t, d_t$ (where $C_t$ is the constant defining $d_t$). This implies that there exists $\alpha >0$ (depending on $t$) such that for any $\epsilon >0$ there exists $C>0$ with
\begin{equation}\label{eq.faff}
\#\{ g \in G: v_{d_\ast} d_\ast(o,g) < T : |td(o,g) + (\theta_{d/d_\ast'}(t) - v_{d_\ast})d_\ast'(o,g)| < \epsilon\} \sim \frac{C e^{\alpha T}}{\sqrt{T}}
\end{equation}
as $T\to\infty$. 

Now note that since $\theta_{d/d_\ast}$ is $C^1$ by Theorem \ref{thm.Man} (and equal to $v_{d_\ast}$ at $t=0$) the function $f: \R \to \R$ defined by
\[
t \mapsto f(t): = \frac{v_{d_\ast} - \theta_{d/d_\ast'}(t)}{t}
\]
is continuous. Also, recall that
\[
\lim_{T\to\infty} f(t) = \sup_{[g] \in \conj(G)} \frac{\ell_{d}[g]}{\ell_{d_\ast'}[g]} \ \ \text{ and } \ \ \lim_{T\to\infty} f(t) = \inf_{[g] \in \conj(G)} \frac{\ell_{d}[g]}{\ell_{d_\ast'}[g]} 
\]
by \cite[Corollary 3.3]{cantrell.tanaka.1}. Combining these two observations implies that there exists $t_0 \in \R$ with $f(t_0) = 1$. When $t = t_0$, equation $(\ref{eq.faff})$ simplifies to the statement in the proposition. This concludes the proof.
\end{proof}
We can now use this lemma to prove the following proposition.

\begin{proposition}\label{prop.decay}
Let $\gfr,\gfr_1,\gfr_2$ be negatively curved, smooth Riemannian metrics on the closed manifold $M$. Suppose that conditions $(1)$ and $(3)$ hold.
Then for any $\epsilon > 0$,
\begin{equation}
\frac{1}{\#\{g \in G: d_\gfr(o,g) < T \}}\#\{ g \in G: d_\gfr(o,g) < T, |d_{\gfr_1}(o,g) - d_{\gfr_2}(o,g)| < \epsilon \} = O\left( \frac{1}{\sqrt{T}}\right)
\end{equation}
as $T\to\infty$.
\end{proposition}
\begin{proof}
Let $s_0 \in \R$ be the abscissa of convergence of the sum
\[
\sum_{g \in G} e^{-sd_\gfr(o,g) - (d_{\gfr_1}(o,g) - d_{\gfr_2}(o,g))}
\]
as $s$ varies. Now, given any $\epsilon >0$, we can assume that there are infinitely many $g\in G$ with $|d_{\gfr_1}(o,g) - d_{\gfr_2}(o,g)| < \epsilon$ otherwise the conclusion of the proposition is trivially true. This assumption implies that $s_0 > 0$.

Again, by Lemma \ref{lem.combo} there exists $C >0$ such that
\[
d(g,h) = 
     \begin{cases}
       s_0d_\gfr(g,h) + (d_{\gfr_1}(g,h) - d_{\gfr_2}(g,h)) + C  &\text{ if $g \neq h$}\\
       0 &\text{ otherwise}
     \end{cases}
\]
belongs to $\Dc_G$. We also define
\[
d_\ast(g,h) = 
     \begin{cases}
       s_0d_\gfr(g,h) + C  &\text{ if $g \neq h$}\\
       0 &\text{ otherwise}
     \end{cases}
\]
where $C >0$ is the same constant used to define $d$. 

Now note that by construction, $d,d_\ast$ have H\"older continuous Busemannn cocycles. Furthermore, (2) holds by Corollary \ref{cor.imp} and by Lemma \ref{lem.nonarith} (since we are assuming $(1)$ and $(3)$) the only $a,b \in \R$ such that $a\ell_d + b \ell_{d_\ast}$ takes values in $\Z$ are $a=b=0$.

Hence we can apply Proposition \ref{prop.asymp} to the pair $d,d_\ast$ to deduce that for any $\epsilon >0$
\begin{equation}\label{eq.candec}
\frac{1}{\#\{g \in G: d_\ast(o,g) < T \}}\#\{ x \in G: d_\ast(o,g) < T, |d(o,g) - d_\ast(o,g)| < \epsilon \} = O\left( \frac{1}{\sqrt{T}}\right)
\end{equation}
as $T\to\infty$. Note however that $d_\ast(o,g) - d(o,g) = d_{\gfr_1}(g,h) - d_{\gfr_2}(g,h)$ and $d_\ast(o,g) = s_0d_\gfr(o,g) +C$ for all $g \neq o$. Hence equation $(\ref{eq.candec})$ is precisely what we wanted to show in the statement of the proposition.
\end{proof}

We are now ready to prove the main theorem of this section.

\begin{proof}[Proof of Theorem \ref{thm.manifolds}]
Take metrics $\gfr_1, \gfr_2$ with distinct marked length spectra, i.e. assume that $(3)$ holds. We will show that
\[
\limsup_{T\to\infty} \frac{\sqrt{T} \, \, \#\{ [g] \in \conj(G): \ell_{\mathfrak{g}}[g] < T, \ell_{\gfr_1}[g] = \ell_{\gfr_2}[g]\}}{\#\{ [g] \in \conj(G): \ell_{\mathfrak{g}}[g] < T \}} < \infty.
\]
As in the proof of Theorem \ref{thm.outermain} we consider two cases depending on the value of $\Lambda_\gfr(\gfr_1,\gfr_2)$.\\

\textit{Case $1$: When $\Lambda_\gfr(\gfr_1,\gfr_2) \neq 0$, i.e. $(1)$ doesn't hold.}
Let $\Cc_\gfr(T) = \{ [g] \in \conj(G): \ell_\gfr[g] < T\}$. By Proposition \ref{prop.ldt} for any $\epsilon >0$  the quantity
\[
\frac{1}{\#\mathcal{C}_\gfr(T)} \#\left\{[g] \in \mathcal{C}_\gfr(T) : \left| \frac{\ell_{\gfr_1}[g] - \ell_{\gfr_2}[g]}{\ell_\gfr[g]} - \Lambda(\gfr_1, \gfr_2) \right| > \epsilon \right\}
\]
decays to $0$ exponentially quickly as $T\to\infty$.
It follows that if $\ell_{\gfr_1}[g] = \ell_{\gfr_2}[g]$ on a subset of conjugacy classes $E$, then the density of this collection with respect to $\ell_\gfr$ must decay to $0$ exponentially quickly.\\

\textit{Case $2$: When $\Lambda_\gfr(\gfr_1, \gfr_2) = 0$, i.e. $(1)$ holds.} 

In this case, Corollary \ref{cor.imp} immediately implies that $(2)$ holds. Hence $(1)$, $(2)$ and $(3)$ hold and 
so we can apply Proposition \ref{prop.decay} to deduce that for any $\epsilon > 0$
\begin{equation}\label{eq.decaysurface}
\frac{1}{\#\{g \in G: d_\gfr(o,g) < T \}}\#\{ g \in G: d_\gfr(o,g) < T, |d_{\gfr_1}(o,g) - d_{\gfr_2}(o,g)| < \epsilon \} = O\left( \frac{1}{\sqrt{T}}\right)
\end{equation}
as $T\to\infty$.

We now consider
\begin{equation*}
\#\{ [g] \in \conj(G): \ell_\gfr[g] < T, \ell_{\gfr_1}[g] =\ell_{\gfr_2}[g]\}
\end{equation*}
which we want to show decays at most like $T^{-1/2} \#\{[g] \in \conj(G): \ell_\gfr[g] <T\}$ as $T\to\infty$.

It suffices to show that
\begin{equation}\label{eq.renorm}
\#\{ [g] \in \conj(G):\text{ $[g]$ is primitive and }  T/2<\ell_\gfr[g] < T, \ell_{\gfr_1}[g] =\ell_{\gfr_2}[g]\} 
\end{equation}
decays at most like $T^{-1/2} \#\{[g] \in \conj(G): \ell_\gfr[g] <T\}$ as $T\to\infty$ since the additional restrictions on $[g]$ cause an exponentially small error term.

As in the proof of Theorem \ref{thm.outermain} we can apply Lemma \ref{lem.hyp} (this time to the pseudometrics  $d_{\gfr_1}, d_{\gfr_2}, d_\gfr$) to deduce that the left hand side of $(\ref{eq.renorm})$ is bounded above by
\begin{equation}\label{eq.bloo}
(M_1T/2 - M_2)^{-1} \#\{g \in G: d_\gfr(o,g) < T+C, |d_{\gfr_1}(o,g) - d_{\gfr_2}(o,g)| < \epsilon + C\}.
\end{equation}
We now note that
\[
\#\{ x \in G: d_\gfr(o,g) < T\} \le T \#\{[g] \in \conj(G): \ell_\gfr[g] < T\}
\]
for all $T >0$ sufficiently large \cite{CoornaertKnieper} and so $(\ref{eq.bloo})$ is bounded above by the quantity in $(\ref{eq.decaysurface})$. This establishes the required estimate for $(\ref{eq.renorm})$ and concludes the proof of the theorem.
\end{proof}


\bibliographystyle{alpha}
\bibliography{Optimal}

\end{document}